\documentclass[twoside,notitlepage,11pt]{article}

\usepackage{amssymb}
\usepackage[leqno]{amsmath}
\usepackage{amsfonts}
\usepackage{amsopn}
\usepackage{amstext}
\usepackage{amsthm}

\usepackage[all]{xy}

\providecommand{\cal}{\mathcal}
\renewcommand{\Bbb}{\mathbb}

\newenvironment{pf}{\begin{proof}}{\end{proof}}

\newcommand{\Aaa}{{\cal{A}}}
\newcommand{\Bee}{{\cal{B}}}

\newcommand{\Ef}{{\cal{F}}}

\newcommand{\Yu}{{\cal{U}}}
\newcommand{\Vee}{{\cal{V}}}

\newcommand{\Qyu}{{\Bbb{Q}}}

\newcommand{\lam}{{\lambda}}
\newcommand{\al}{\alpha}

\newcommand{\sig}{\sigma}

\renewcommand{\phi}{\varphi}
\renewcommand{\rho}{\varrho}

\newcommand{\rest}{\restriction}

\newcommand{\ntr}{n\in\omega}

\newcommand{\loe}{\leqslant}
\newcommand{\goe}{\geqslant}

\newcommand{\subs}{\subseteq}

\newcommand{\nnempty}{\ne\emptyset}

\newcommand{\cl}{\operatorname{cl}}

\newcommand{\w}{\operatorname{w}}

\newcommand{\id}{\operatorname{id}}
\newcommand{\cf}{\operatorname{cf}}

\newcommand{\supp}{\operatorname{suppt}}

\newcommand{\liminv}{\varprojlim}

\newcommand{\concat}{{}^\smallfrown}

\newtheorem{tw}{Theorem}[section]
\newtheorem{wn}[tw]{Corollary}
\newtheorem{lm}[tw]{Lemma}
\newtheorem{prop}[tw]{Proposition}

\theoremstyle{definition}

\newtheorem{ex}[tw]{Example}
\theoremstyle{remark}

\newcommand{\setof}[2]{\{#1\colon #2\}}

\newcommand{\seqof}[2]{\langle #1\colon #2\rangle}
\newcommand{\sett}[2]{\{#1\}_{#2}}
\newcommand{\sn}[1]{\{#1\}} 
\newcommand{\dn}[2]{\{#1,#2\}} 
\newcommand{\pair}[2]{\langle #1, #2 \rangle} 
\newcommand{\map}[3]{#1\colon #2 \to #3} 
\newcommand{\img}[2]{#1[#2]} 
\newcommand{\inv}[2]{{#1}^{-1}[#2]} 

\newcommand{\dpower}[2]{[#1]^{#2}}

\newcommand{\rloe}{\preceq}

\newcommand{\suppt}{\supp}

\newcommand{\Sig}{\Sigma}

\newcommand{\II}{\ensuremath{\mathbb I}}

\newcommand{\iso}{\cong}

\newcommand{\invsys}[5]{\langle {#1}_{#4};{#2}_{#4}^{#5};#3 \rangle}
\newcommand{\rinv}[4]{\langle{{#1}_{#3}^{#4}};{#2}\rangle} 
\renewcommand{\S}{\mathbb S}
\newcommand{\diag}{\Delta} 
\newcommand{\Diag}{{\Delta}} 
\newcommand{\cube}[1]{\ensuremath{[0,1]^{#1}}} 
\newcommand{\til}{\tilde}

\title{Small Valdivia compact spaces}
\author{
{\sc Wies{\l}aw Kubi\'s}\\ \\
University of Silesia\\
Katowice, Poland\\
\texttt{kubis@ux2.math.us.edu.pl}
\and
{\sc Henryk Michalewski}\\ \\
Warsaw University\\
Warsaw, Poland\\
\texttt{henrykm@duch.mimuw.edu.pl}
}

\begin{document}
\maketitle

\begin{abstract} We prove a preservation theorem for the class of Valdivia compact spaces, which involves inverse sequences of ``simple'' retractions. Consequently, a compact space of weight $\loe\aleph_1$ is Valdivia compact iff it is the limit of an inverse sequence of metric compacta whose bonding maps are retractions.
As a corollary, we show that the class of Valdivia compacta of weight $\loe\aleph_1$ is preserved both under retractions and under open $0$-dimensional images. Finally, we characterize the class of all Valdivia compacta in the language of category theory, which implies that this class is preserved under all continuous weight preserving functors.

\noindent{\bf AMS Subject Classification (2000)} Primary: 54D30; Secondary: 54C15, 54B35.

\noindent{\bf Keywords and phrases:} Valdivia compact space, retraction, $\Sigma$-product, inverse system, r-skeleton.
\end{abstract}

\section{Introduction}

A topological space is called {\em Valdivia compact\/}, if it is homeomorphic to the closure in a Tikhonov cube of a subset of a $\Sigma$-product. This is a natural generalization of {\em Corson compact spaces\/}, which are defined as compact subsets of $\Sigma$-products.
The class of Valdivia compact spaces was introduced by Argyros, Mercourakis and Negrepontis in \cite{AMN} and further investigated by Valdivia \cite{Valdivia, Valdivia3}. The name {\em Valdivia compact} was introduced by Deville and Godefroy in \cite{DG}. 
Valdivia compacta were extensively studied by Kalenda and we refer to his article \cite{Kalenda} for a survey of results. It is well known that the class of Corson compacta is stable under continuous images. On the other hand, there are simple examples of continuous images of Valdivia compacta which are not Valdivia compact. An interesting result of Kalenda \cite{Kalenda4} says that every Valdivia compact space which is not Corson compact, has a two-to-one continuous map onto a non-Valdivia compact space. 

The main subject of this paper is to study Valdivia compact spaces of weight $\loe\aleph_1$ (we call these spaces {\em small Valdivia compacta}). A typical phenomenon occurs when studying various classes of nonmetrizable compact spaces: namely, a class restricted to spaces of weight $\loe\aleph_1$ has a very simple structure and above $\aleph_1$ the structure becomes much more complicated. 
In this context, we prove some positive results about the stability of the class of small Valdivia compacta. It turns out that small Valdivia compacta are precisely those spaces which can be obtained as limits of inverse sequences of metric compacta whose bonding maps are retractions. 
In fact, we prove a more general preservation result which involves inverse sequences of {\em simple} retractions (the definition is given in Section \ref{prosteretrakcje}).
We apply this result to give positive answers to two questions on the stability of Valdivia compacta, in the case of spaces of weight $\loe\aleph_1$. The mentioned preservation result shows also that retracts of Tikhonov or Cantor cubes are Valdivia compact. A stronger result has recently been obtained by Leiderman and the first author in \cite{KL}, namely retracts of products of metric compacta belong to the class of {\em semi-Eberlein spaces} and the latter class is properly contained in the class of Valdivia compacta.

A question of Kalenda \cite{Kalenda} asks whether every open continuous image of a Valdivia compactum is Valdivia compact. He gives an affirmative answer in the case where the image has a dense set of $G_\delta$ points, see \cite{Kalenda5, Kalenda}. 
Uspenskij and the first author described in \cite{KU} an example of a compact connected Abelian group of weight $\aleph_1$ which is not Valdivia compact. Such a group is an open (epimorphic) image of a Valdivia compactum (a product of compact metric spaces), which gives a negative answer to Kalenda's question. On the other hand, in this paper we show that a small open $0$-dimensional image of a Valdivia compact space is Valdivia compact.
Another question from \cite{Kalenda} asks whether $X$ is Valdivia compact provided $X\times Y$ is Valdivia compact for some space $Y$. We give a partial positive answer to a more general question whether the class of Valdivia compacta is closed under retractions. Namely, we show that a retract of a Valdivia compactum is Valdivia compact provided its weight is $\loe\aleph_1$. 

The last part of the paper explores the properties of inverse systems of retractions. As a result, we give a characterization of the class of all Valdivia compacta in terms of certain inverse systems of retractions. This characterization belongs to the language of category theory and therefore as an application we obtain that every continuous weight preserving functor on the class of compact spaces preserves Valdivia compacta.

\section{Preliminaries}

\subsection{Notation and basic definitions}
All topological spaces are assumed to be Hausdorff. By a ``map" we mean a continuous map, unless otherwise indicated. The closure of a set $A$ is denoted by $\cl A$ or $\cl_XA$, if it is not clear from the context, which topological space we have in mind.
We say that $A$ is {\em countably closed} if $\cl M\subs A$ whenever $M\subs A$ is countable. A subset $D$ of a topological space $X$ is {\em $\kappa$-monolithic} if for every $A\in\dpower D{\loe\kappa}$, $\cl_XA$ is contained in $D$ and its net-weight is $\loe\kappa$. As we shall consider $\kappa$-monolithic sets in compact spaces, ``net-weight" can be replaced by ``weight" in the above definition.

A Tikhonov cube of weight $\kappa$ is denoted by $\cube S$, where $S$ is any set of cardinality $\kappa$. If $T\subs S$ then we denote by $x\mid T$ the element $(x\rest T)\concat 0_{S\setminus T}\in\cube S$, where $0_{S\setminus T}$ is the constant $0$ function in $\cube{S\setminus T}$ and $\concat$ denotes the concatenation of functions. We also write $X\mid T$ for the set $\setof{x\mid T}{x\in X}$, where $X\subs \cube S$.
The map $\map{j_T^S}{\cube T}{\cube S}$ defined by $j_T^S(x)=x\concat 0_{S\setminus T}$ is called the {\em canonical embedding} of $\cube T$ into $\cube S$. We shall often identify $\cube T$ with a subspace of $\cube S$, meaning that canonical embedding. The standard projection (i.e. the restriction map) from $\cube S$ onto $\cube T$ is denoted by $\pi^S_T$. 

Given two maps $\map fXY$ and $\map gXZ$, the {\em diagonal} of $f,g$ is the map $\map{f\diag g}X{Y\times Z}$ defined by $(f\diag g)(x)=\pair{f(x)}{g(x)}$.
The disjoint topological sum of spaces $X,Y$ is denoted by $X\oplus Y$. 
Ordinals will be denoted by Greek letters and treated as sets of their predecessors linearly ordered by the relation $\in$. Cardinals are particular examples of ordinals, although we usually denote them by using the $\aleph$-notation. Given an ordinal $\delta$, the set $\delta+1$ is a compact linearly ordered space. In particular, $\omega+1$ (which is denoted also by $\aleph_0+1$) is a convergent sequence together with its limit. The cardinality of a set $A$ is denoted by $|A|$.

A family $\Aaa$ of subsets of a space $X$ is {\em $T_0$ separating} if $x,y\in X$ and $x\ne y$ imply that there is $A\in\Aaa$ such that $|\dn xy\cap A|=1$. Every $T_0$ separating family $\Yu$ of open $F_\sig$ subsets of a compact space $X$ induces an embedding $\map fX{\cube\Yu}$ defined by $f(x)(u)=h_u(x)$, where $\map{h_u}X{[0,1]}$ is a fixed continuous function such that $u=\inv{h_u}{(0,1]}$, $u\in\Yu$.

A {\em retraction} is, by definition, a map $\map fXY$ which has a {\em right inverse} $\map gYX$, i.e. $fg=\id_Y$. In this case $f$ is a surjection, $g$ is an embedding and $\map{gf}XX$ is a self-map which is identity on its range. Conversely, if $X$ is compact and $\map rXX$ is such that $r^2=r$ then $\map rX{\img rX}$ is a retraction in the above sense. We then say that $r$ is an {\em internal retraction}. We shall use the following well known fact: if $\map fXY$ is an open surjection, $X$ is a complete metric space and $Y$ is a $0$-dimensional compact space then $f$ is a retraction. The proof is a direct application of Michael's selection theorem for the multivalued map $f^{-1}$, which is lower semi-continuous in this case.

We shall use the method of elementary substructures, in order to avoid tedious ``closing-off" constructions. We refer to \cite{Dow} for a survey on the use of this method in set theoretic topology.
Given an uncountable cardinal $\chi$ we denote by $H(\chi)$ the class of all sets which are hereditarily of cardinality $<\chi$. That is $x\in H(\chi)$ iff $|x|<\chi$, $|y|<\chi$ for every $y\in x$, $|z|<\chi$ for every $z\in y\in x$ and so on. It can be shown easily that $H(\chi)$ is a set and for every formula $\phi(x_1,\dots,x_n)$ (with parameters $x_1,\dots,x_n$) which is satisfied in the universe of set theory we are working in, there exists a cardinal $\chi$ such that the structure $\pair{H(\chi)}{\in}$ contains the parameters $x_1,\dots,x_n$ and satisfies $\phi(x_1,\dots,x_n)$. This is called the {\em reflection principle}. Given finitely many objects, say $a_1,\dots, a_n$ and having in mind finitely many formulas which speak about these objects, there is always a ``big enough" cardinal $\chi$ such that $a_1,\dots,a_n\in H(\chi)$ and all the formulas are satisfied in $\pair{H(\chi)}{\in}$ provided they are satisfied in the universe. We shall ignore the details, saying ``fix a big enough $\chi$", having in mind a cardinal such that $H(\chi)$ contains all the relevant objects and reflects all the finitely many relevant formulas. We write $M\rloe H(\chi)$ for ``$M$ is an elementary substructure of $\pair{H(\chi)}{\in}$". Every $M\rloe H(\chi)$ has the following properties which we shall use throughout the paper:
\begin{enumerate}
\item[(i)] if $x\in M$ and $|x|\loe\aleph_0$ then $x\subs M$. 
\item[(ii)] if $x_0,\dots, x_{n-1}\in M$ then $\{x_0,\dots,x_{n-1}\}\in M$. 
\end{enumerate}
A particular case of the well known theorem of Skolem-L\"owenheim says that for every infinite set $S\subs H(\chi)$ there exists $M\rloe H(\chi)$ such that $S\subs M$ and $|S|=|M|$.
We shall demonstrate the use of elementary substructures in the next subsection, reproving some properties of $\sig$-complete inverse systems.

\subsection{Inverse systems}
We shall consider inverse systems of compact space whose all bonding mappings are surjections. An inverse system indexed by $\Sigma$ will be denoted by $\S=\invsys Xp\Sigma st$, where $X_s$ is the $s$-th space in the system and $\map{p^t_s}{X_t}{X_s}$ is the bonding map, $s\loe t$. Formally, the limit of $\S$ is a pair $\pair X{\sett{p_s}{s\in\Sigma}}$, where $\map{p_s}X{X_s}$ and for every space $Y$ and for every collection of maps $\setof{f_s}{s\in\Sigma}$ such that $\map{f_s}Y{X_s}$ and $s\loe t\implies f_s= p^t_sf_t$, there exists a unique map $\map gYX$ such that $f_s=p_s g$ holds for every $s\in\Sigma$. Mappings $p_s$ are called {\em projections}. In most cases it is clear from the context what the projections are, then we say that $X$ is the limit of $\S$ and we write $X=\liminv\S$, instead of $\pair X{\sett{p_s}{s\in\Sigma}}=\liminv\S$. If $T$ is a directed subset of $\Sigma$ then $\S\rest T:=\invsys XpTst$ is again an inverse system. $\S$ is {\em continuous} if $X_s=\liminv(\S\rest T)$, whenever $T\subs \Sigma$ is directed and such that $s=\sup T$. It is well known that $\liminv \S=\liminv(\S\rest T)$, whenever $T$ is cofinal in $\Sigma$. We say that $\S=\invsys Xp\Sigma st$ is {\em $\sig$-complete} if every countable directed set $T\subs\Sigma$ has the least upper bound in $\Sigma$ (i.e. $\Sigma$ is {\em $\sig$-complete}) and $X_t=\liminv(\S\rest T)$, where $t=\sup T$. Given a $\sig$-complete directed poset $\Sigma$, we say that $T\subs \Sigma$ is {\em $\sig$-closed} in $\Sigma$ if $\sup M\in T$ for every countable directed set $M\subs T$.

Given two maps $\map fXY$ and $\map gXZ$, we say that $f$ {\em factors through} $g$ if there exists $\map hZY$ such that $f= hg$. If both $f,g$ are quotient maps then such a map $g$ is unique and it is automatically continuous. 

The following properties of $\sig$-complete inverse systems will be used throughout the paper.

\begin{prop}\label{faktor} Let $\S=\invsys Xp\Sigma st$ be a $\sig$-complete inverse system of metric compacta and let $X=\liminv\S$. Then for every map $\map fXY$ into a second countable space there exists $\delta\in\Sigma$ such that $f$ factors through $p_\delta$, i.e. there exists $\map g{X_\delta}Y$ such that $f=g p_\delta$.
\end{prop}

\begin{pf} Fix a big enough cardinal $\chi$ and a countable $M\rloe H(\chi)$ which knows $\S$ and $f$. By elementarity, $\Sigma\cap M$ is directed. Let $\delta=\sup(\Sigma\cap M)$. We claim that $f$ factors through $p_\delta$. 

Fix $x_0,x_1$ such that $f(x_0)\ne f(x_1)$. We have to show that $p_\delta(x_0)\ne p_\delta(x_1)$. Fix disjoint open sets $u_0, u_1\subs Y$ such that $f(x_i)\in u_i$. We may assume that $u_0,u_1\in M$, because $M$ knows some countable base of $Y$. Now $\inv f{u_0}$ is an open $F_\sig$ subset of $X$ and therefore by compactness it is the union of a countable family of open sets of the form $\inv{p_s}w$, where $s\in\Sig$ and $w\subs X_s$. By elementarity, we may assume that such a family belongs to $M$ and consequently is contained in $M$. In particular, there are $s_0\in\Sigma\cap M$ and an open set $w_0\subs X_{s_0}$ such that $x_0\in \inv {p_{s_0}}{w_0}\subs \inv f{u_0}$. 
By the same argument, there are $s_1\in\Sigma\cap M$ and an open set $w_1\subs X_{s_1}$ such that $x_1\in \inv {p_{s_1}}{w_1}\subs \inv f{u_1}$. Modifying the sets if necessary and using the fact that $\Sigma\cap M$ is directed, we may assume that $s_0=s_1=s\in M$. Then $w_0\cap w_1=\emptyset$ and consequently $p_\delta(x_0), p_\delta(x_1)$ are separated by disjoint sets $\inv{(p^\delta_s)}{w_i}$, $i=0,1$. Thus $p_\delta(x_0)\ne p_\delta(x_1)$.
\end{pf}

\begin{prop}\label{wsegasgpajspjwqprzmv}
Assume $X=\liminv\invsys Xp\Sigma s{s'}$, $Y=\liminv\invsys Yq\Sigma s{s'}$, where both inverse systems are $\sig$-complete and $X_s,Y_s$ are metric compacta for every $s\in\Sigma$. Let $\map fXY$. Then there are a $\sig$-closed cofinal set $T\subs \Sigma$ and a family of maps $\sett{f_t}{t\in T}$ such that $\map{f_t}{X_t}{Y_t}$ and $q_tf = f_tp_t$ holds for every $t\in T$. If moreover the map $f$ is open then we may assume that $f_t$ is open for every $t\in T$.
\end{prop}

\begin{pf} Since the projections are quotient maps (recall that we consider inverse systems whose bonding maps are onto), we do not have to worry about the continuity of the maps $f_t$. Denote by $T$ the set of all $t\in \Sigma$ such that $q_tf$ factors through $p_t$, i.e. such that there is (necessarily uniquely defined) $f_t$ satisfying $q_tf = f_tp_t$ and $f_t$ is open if $f$ is so. Then $T$ is $\sig$-closed, by the continuity of both systems (note that the limit of an inverse system of open maps is open). 
It remains to check that $T$ is cofinal in $\Sigma$. 

Fix a big enough cardinal $\chi$, fix $s\in \Sigma$ and find a countable $M\rloe H(\chi)$ which contains enough information and such that $s\in M$. Let $\delta=\sup(\Sigma\cap M)$. Then $s<\delta$ and we claim that $\delta\in T$. Fix $x_0,x_1\in X$ such that $q_\delta f(x_0)\ne q_\delta f(x_1)$. Then $q_tf(x_0)\ne q_tf(x_1)$ for some $t\in \Sigma\cap M$ and by elementarity we can find $\phi\in C(X_t)\cap M$ such that $\phi q_tf(x_0)<\phi q_tf(x_1)$. By Proposition \ref{faktor}, $\phi q_tf$ factors through $p_r$ for some $r$. Since $\phi q_tf\in M$, we may assume that $r\in M$. Thus $r<\delta$ and $\phi q_tf$ factors also through $p_\delta$, i.e. there exists a map $g$ such that $gp_\delta=\phi p_tf$. Hence $gp_\delta(x_0)<gp_\delta(x_1)$, which shows that $p_\delta(x_0)\ne p_\delta(x_1)$.
It follows that $q_\delta f$ factors through $p_\delta$.

Assume now that $f$ is open and denote by $f_\delta$ the map which realizes the factorization (where $\delta$ is as above). We claim that $f_\delta$ is open. Fix a basic open set $u=\inv{(p^\delta_t)}w$, where $w\subs X_t$ and $t\in \Sigma\cap M$. Since $M$ knows a countable base of $X_t$, we may assume that $w\in M$. Let $v=f \inv{p_t}w=f\inv{p_\delta}u$. Then $v\in M$ is an open $F_\sig$ set and therefore by elementarity there is a countable collection $\Vee$ of open sets of the form $\inv {q_r}a$ whose union is $v$. By elementarity, we may assume that $\Vee\in M$ and therefore $v=\bigcup_{\ntr}\inv{q_{r_n}}{a_n}$, where $\setof{r_n}{\ntr}\subs M$. Let $b_n=\inv{(q^\delta_{r_n})}{a_n}$. Then $v=\inv {q_\delta}b$, where $b=\bigcup_{\ntr}b_n$. Finally, we have $b=q_\delta f \inv{p_\delta}{u}= f_\delta p_\delta \inv{p_\delta}{u} = \img{f_\delta}{u}$, therefore $\img{f_\delta}u$ is open. 
\end{pf}

\begin{prop}[Spectral Theorem, cf. \cite{Szczepin}]\label{asfjsifjioqwrqr} 
Let $X=\liminv\invsys Xp\Sigma s{s'}$, $Y=\liminv\invsys Yq\Sigma s{s'}$, where both inverse systems are $\sig$-complete and $X_s, Y_s$ are metric compact for every $s\in\Sigma$.
Assume further that $\map fXY$ is a retraction. Then there are a $\sig$-closed cofinal set $T\subs\Sigma$ and a family of maps $\sett{f_t}{t\in T}$ such that $q_tf = f_tp_t$ and $f_t$ is a retraction for every $t\in T$.
\end{prop}

\begin{pf}
Assume first that $f$ is a retraction and let $\map iYX$ be its right inverse. By Proposition \ref{wsegasgpajspjwqprzmv} applied to both $f$ and $i$, there are a $\sig$-closed cofinal set $T\subs \Sigma$ and families of maps $\sett{f_t}{t\in T}$, $\sett{i_t}{t\in T}$ such that $q_tf = f_tp_t$ and $p_ti = i_tq_t$ for every $t\in T$. Then, given $t\in T$, we have $f_ti_tq_t = f_tp_ti = q_tfi = q_t$ and therefore $f_ti_t=\id_{Y_t}$, because $q_t$ is a surjection.
\end{pf}

\subsection{$\Sigma$-products}

Given a set $\kappa$ (not necessarily a cardinal), we denote by $\Sigma(\kappa)$ the $\Sigma$-product modelled on $\kappa$, which by definition consists of all $x\in\cube\kappa$ such that $\suppt(x):=\setof{\al\in \kappa}{x(\al)\ne0}$ is countable. Clearly, $\Sigma(\kappa)$ is $\aleph_0$-monolithic in $\cube\kappa$.
It is well known that for every set $\kappa$, $\Sigma(\kappa)$ has countable tightness. Another result, due to Corson \cite{Corson}, says that every two disjoint relatively closed sets $A,B\subs\Sigma(\kappa)$ have disjoint closures in $\cube\kappa$. This implies that $\Sigma(\kappa)$ is a normal space and for every relatively closed set $S\subs\Sigma(\kappa)$ the space $\cl_{\cube\kappa} S$ is (homeomorphic to) the \v{C}ech-Stone compactification of $S$. Below we state and prove a more general property of $\Sigma$-products which implies the above result and which will be used later.

\begin{lm}\label{wiewior1} Let $X\subs\cube\kappa$ be such that $X\cap \Sigma(\kappa)$ is dense in $X$ and let $\map fXY$ be a continuous map. Let $\Bee$ be an open base for $\img fX$. Assume $\chi$ is big enough and $M\rloe H(\chi)$ is such that $f,\kappa\in M$ and $\Bee\subs M$. Let $S=M\cap \kappa$. Then
\begin{enumerate}
	\item[$(a)$] $x\mid S\in \cl_{\cube\kappa} X$ for every $x\in X$.
	\item[$(b)$] If\/ $x\in X$ and $x\mid S\in X$ then $f(x)=f(x\mid S)$.
\end{enumerate}
\end{lm}

\begin{pf} (a) Fix $x\in X$ and fix a basic neighborhood $u$ of $x\mid S$. Then $u$ is defined by a finite function, say $\psi$, whose values are open rational intervals. Let $u'$ denote the neighborhood defined by the function $\psi$ restricted to $S$. Then $u'\in M$.
Now $M\models u'\cap X\nnempty$ so there exists $y\in M$ with $y\in \Sigma(\kappa)\cap u'\cap X$. Thus $\supp(y)\subs M$ (because it is countable) and therefore $\supp(y)\subs S$. It follows that $y\in u$, so $u\cap X\nnempty$. As $u$ is arbitrary, we get $x\mid S\in\cl X$.

(b) Suppose $x,x\mid S\in X$ and $f(x)\ne f(x\mid S)$. Let $v,w\in\Bee$ be disjoint such that $f(x)\in v$ and $f(x\mid S)\in w$. There exist basic neighborhoods $u_1,u_2$ of $x$ and $x\mid S$ respectively, such that $\img f{u_1}\subs v$ and $\img f{u_2}\subs w$.
Assuming $u_i$ is defined by a finite function $\psi_i$ whose values are open rational intervals and, shrinking both neighborhoods if necessary, we may also assume that $\psi=\psi_1\rest S=\psi_2\rest S$. Let $u$ denote the basic open set defined by $\psi$. Then $u\in M$.
Furthermore, $M\models \inv f v\cap u\nnempty$. Thus, there exists $y\in M$ with $y\in \Sigma(\kappa)\cap X\cap u$ such that $f(y)\in v$. But $\suppt(y)\subs M$; hence $y(i)=0$ for $i\notin S$ and therefore $y\in u_2$, because $(x\mid S)(i)=0\in \psi_2(i)$ for $i\notin S$.
It follows that $f(y)\in w$, a contradiction.
\end{pf}

\subsection{Valdivia compacta}
A space $X$ is said to be {\em Valdivia compact} \cite{DG} if for some set $T$ there exists an embedding $\map hX{\cube T}$ such that 
$$\img hX=\cl_{\cube T}\Bigl(\Sigma(T)\cap\img hX\Bigr).$$
Such a map $h$ will be called a {\em Valdivia embedding}. More generally, we say that $\map fX{\cube T}$ is a {\em Valdivia map} if $\Sigma(T)\cap \img fX$ is dense in $\img fX$ (so $f$ does not have to be one-to-one).

Clearly, every Corson compact space is Valdivia compact. In fact, Corson compacta are precisely those Valdivia compact spaces which have countable tightness (because $\Sigma$-products are countably closed). Below we recall some basic properties of Valdivia compacta. For completeness, we give the proofs.

\begin{prop}\label{koty} Assume $X$ is Valdivia compact and $\kappa=\w(X)$. Then there exists a Valdivia embedding of $X$ into $\cube\kappa$. \end{prop}

\begin{pf} Assume $X\subs\cube\lam$ and the inclusion is a Valdivia embedding.
Let $D=\setof{d_\al}{\al<\kappa}\subs \Sigma(\lam)\cap X$ be dense in $X$. 
Let $T=\bigcup_{\al<\kappa}\suppt(d_\al)$. Then $|T|\loe\kappa$ and, identifying $\cube T$ with a subset of $\cube \lam$, we have $X\subs\cube T$.
\end{pf}

\begin{prop}[\cite{AMN}] Assume $X$ is Valdivia compact of weight $\kappa$. Then $X=\liminv\S$, where $\S=\invsys Xr\kappa\al\beta$ is a continuous inverse sequence such that each $X_\al$ is Valdivia compact of weight $\loe\al+\aleph_0$ and each $r^\beta_\al$ is a retraction. \end{prop}

\begin{pf} Assume $X\subs \cube\kappa$ and the inclusion is a Valdivia embedding. Using Lemma \ref{wiewior1}, find a continuous increasing sequence $\seqof{S_\al}{\al<\kappa}$ of subsets of $\kappa$ such that $|S_\al|\loe\al+\aleph_0$ and $x\mid S_\al\in X$ for every $x\in X$, $\al<\kappa$. Set $X_\al=X\mid S_\al$ and $r^\beta_\al=R_{S_\al}\rest X_\beta$, where $R_{S_\al}(x)=x\mid S_\al$. Then $\S=\invsys Xr\kappa\al\beta$ is as required and $X=\liminv \S$.
\end{pf}

\begin{prop}[\cite{Gl}]\label{sigmas} Let $X\subs \cube\kappa$ be a Valdivia inclusion and let $\Sigma=\Sigma(\kappa)\cap X$. Then $X=\beta\Sigma$.
\end{prop}

\begin{pf} We have to show that every two disjoint relatively closed sets $A,B\subs\Sigma$ have disjoint closures in $X$. Suppose $p\in\cl_XA\cap \cl_XB$ and fix a countable $M\rloe H(\chi)$ such that $A,B,p \in M$. Let $S=\kappa\cap M$. Lemma \ref{wiewior1}(a) applied both to $A\cup\sn p$ and $B\cup\sn p$ says that 
$$p\mid S\in\cl(A\cup\sn p)\cap \cl(B\cup\sn p)=(\cl A\cap \cl B)\cup\sn p,$$
where $\cl$ denotes the closure in $\cube\kappa$.
Thus $p\mid S\in A\cap B$, because $S$ is countable; a contradiction. \end{pf}

Note that the set $\Sigma$ from the above proposition is $\aleph_0$-monolithic. The existence of a dense $\aleph_0$-monolithic subset is a property preserved by continuous images. Thus, there are obvious examples of spaces which are not continuous images of any Valdivia compactum, like e.g. non-metrizable compactifications of $\omega$.

\section{Inverse systems with right inverses I}

In this section we study inverse systems whose bonding mappings are retractions.

Let $\S=\invsys Xr\Sigma st$ be an inverse system. A collection of embedddings $\setof{i_s^t}{s\loe t,\;s,t\in\Sigma}$, where $\map{i_s^t}{X_s}{X_t}$, satisfying $i_s^s=\id_{X_s}$, $r_s^t i_s^t=\id_{X_s}$ and $i^r_t i^t_s=i^r_s$ whenever $s\loe t\loe r$, will be called a {\em right inverse of\/} $\S$. We shall write briefly $\rinv i\Sigma st$ instead of $\setof{i_s^t}{s\loe t,\;s,t\in\Sigma}$. Clearly, if $\S$ has a right inverse then all bonding maps in $\S$ are retractions. In fact, all projections from the limit are also retractions, by the following lemma.

\begin{lm}\label{elm} Assume $X=\liminv\S$, where $\S=\invsys Xr\Sigma st$ is an inverse system with a right inverse $\rinv i\Sigma st$. There are uniquely determined embeddings $\map{i_s}{X_s}X$ such that $r_s i_s=\id_{X_s}$ and $i_t i_s^t=i_s$ for every $s\loe t$ in $\Sigma$. \end{lm}

\begin{pf} Fix $s\in \Sigma$. Take the usual representation of $X$ as a subset of the product $\prod_{t\goe s}X_t$ and define
$$i_s(x)(t)=i^t_s(x),\quad x\in X_s,\; s\loe t.$$
Given $s<t<u$, we have $r^u_t(i_s(x)(u))=r^u_t i^u_t i^t_s(x)=i_s(x)(t)$,
so $i_s$ is a well defined map $\map{i_s}{X_s}X$. Clearly, $i_s$ is continuous, because $r_t i_s = i^t_s$ is continuous for every $t\goe s$. Moreover $r_s i_s=\id_{X_s}$. Finally, if $s<t$ and $u\goe t$ then
$$(i_t i^t_s(x))(u) = i^u_t(i^t_s(x)) = i^u_s(x) = i_s(x)(u).$$
Thus $i_t i^t_s = i_s$.
\end{pf}

Given an inverse system $\invsys Xr\Sigma st$ with a right inverse $\rinv i\Sigma st$, we shall use the induced embeddings $i_s$, without refering to the above lemma explicitly.

\begin{lm}\label{prawy}
Let $\S=\invsys Xr\kappa\al\beta$ be a continuous inverse sequence such that each $r^{\al+1}_\al$ is a retraction. Then $\S$ has a right inverse. 
\end{lm}

\begin{pf}
We construct a right inverse $\rinv i\kappa\al\beta$ using induction on $\delta<\kappa$. Assume that $i^\beta_\al$ have already been defined for $\al\loe\beta<\delta$. 

If $\delta$ is a limit ordinal then we use Lemma \ref{elm} and the continuity of the sequence.
So suppose that $\delta=\rho+1$. As $r^{\rho+1}_\rho$ is a retraction, choose $\map i{X_\rho}{X_{\rho+1}}$ which is a right inverse of $r^{\rho+1}_\rho$ and define
$i^{\rho+1}_\al=i i^\rho_\al$ for $\al<\rho+1$. In particular, $i^{\rho+1}_\rho=i$.
We have
$$r^{\rho+1}_\al i^{\rho+1}_\al = r^{\rho+1}_\al i i^\rho_\al = r^\rho_\al r^{\rho+1}_\rho i i^\rho_\al = r^\rho_\al i^\rho_\al=\id_{X_\al}$$
and
$$i^{\rho+1}_\beta i^\beta_\al=i i^\rho_\beta i^\beta_\al = i i^\rho_\al = i^{\rho+1}_\al,$$
whenever $\al<\beta<\rho+1$. This completes the proof.
\end{pf}

In the context of Lemmas \ref{elm} and \ref{prawy}, let us mention that Shapiro \cite{Shapiro} constructed (assuming $2^{\aleph_0}=\aleph_1$ and $2^{\aleph_1}=\aleph_2$) an example of a $0$-dimensional $\kappa$-metrizable space of weight $\aleph_2$, which has no nontrivial convergent sequences. Such a space is the limit of a $\sig$-complete inverse system of metric compact spaces whose all bonding mappings are open, therefore retractions. However, the projections from the limit are open mappings which are not retractions.

\begin{lm}\label{brzim} Let $\S=\invsys Xr\Sigma st$ be an inverse system with a right inverse $\rinv i\Sigma st$ and let $X=\liminv\S$, $R_s=i_s r_s$. Then
\begin{enumerate}
	\item $s\loe t\implies R_s R_t=R_s=R_t R_s$.
	\item $x=\lim_{s\in\Sigma}R_s(x)$ for every $x\in X$.
\end{enumerate}
\end{lm}

\begin{pf} Property 1 is clear. Fix a neighborhood $u$ of $x$. We can assume that $u=\inv{(r_s)}v$ for some $s\in\Sigma$, and an open set $v\subs X_s$. Fix $t\goe s$. Then
$$r_s(R_t(x))=r_s(i_t r_t(x))= r^t_s r_t i_t r_t(x)=r^t_s r_t(x)=r_s(x)\in v.$$
Thus $R_t(x)\in u$ for every $t\goe s$. \end{pf}

\begin{lm}\label{antybrzim} Assume $\setof{R_s}{s\in\Sig}$ is a family of internal retractions of a compact space $X$ such that $\Sig$ is a directed partially ordered set and conditions 1, 2 of Lemma \ref{brzim} hold. Then $X=\liminv\invsys XR\Sig st$, where $X_s=\img{R_s}X$ and $R^t_s=R_s\rest X_t$.
\end{lm}

\begin{pf} If $r<s<t$ then 
$$R^s_r R^t_s = (R_r\rest X_s) (R_s\rest X_t) = (R_r R_s)\rest X_t = R_r \rest X_t = R^t_r,$$
thus $\S=\invsys XR\Sig st$ is indeed an inverse system and if $s<t$ then $\img{R^t_s}{X_t}=X_s$, because $X_s\subs X_t$ by condition 1. Let $\map{\til R_s}{\liminv\S}{X_s}$ denote the projection. The collection $\setof{R_s}{s\in\Sig}$ is compatible with $\S$ which determines a continuous map $\map fX{\liminv\S}$ such that $\til R_s f = R_s$ for every $s\in\Sig$. Condition 2 implies that $f$ is one-to-one.
It remains to show that $\img fX = \liminv\S$. Fix $y\in\liminv\S$ and pick
$$x\in\bigcap_{t\in\Sig}\cl\setof{\til R_s(y)}{s\goe t}.$$
The above set is nonempty by the compactness of $X$ and by the directedness of $\Sig$.
Observe that if $t<s$ then $R_t \til R_s(y)= R^s_t \til R_s(y)=\til R_t(y)$, whence
$$R_t(x)\in \cl\img{R_t}{\setof{\til R_s(y)}{s\goe t}}=\sn{\til R_t(y)}.$$
Thus $\til R_t f(x)=R_t(x)=\til R_t(y)$, which means that $f(x)=y$.
\end{pf}

\section{Simple retractions}\label{prosteretrakcje}

We define a notion of a simple retraction and we show that the class of Valdivia compact spaces is stable under limits of inverse sequences with simple retractions. 
As a consequence, we obtain a characterization of Valdivia compacta of weight $\aleph_1$, which allows to answer the question on retractions and open images in case where the image space has weight $\aleph_1$. Another corollary says that retracts of Cantor/Tikhonov cubes are Valdivia compact.

\begin{lm}\label{smok} Assume $X=\liminv\S$, where $\S=\invsys Xr\kappa\al\beta$ is a continuous inverse sequence of compact spaces with a right inverse $\rinv i\kappa\al\beta$.
Assume further that there is a collection of Valdivia embeddings $\sett{\map{h_\al}{X_\al}{\cube{S_\al}}}{\al<\kappa}$ such that for every $\al<\beta<\kappa$ the following diagrams commute:
$$\xymatrix{{X_\beta} \ar[r]^-{h_\beta}\ar[d]_{r^\beta_\al} & {\cube{S_\beta}}\ar[d]^{\pi^{S_\beta}_{S_\al}}\\
{X_\al} \ar[r]^-{h_\al} & {\cube{S_\al}}
} \qquad\qquad
\xymatrix{{X_\beta} \ar[r]^-{h_\beta} & {\cube{S_\beta}}\\
{X_\al} \ar[u]^{i^\beta_\al} \ar[r]^-{h_\al} & {\cube{S_\al}} \ar[u]_{j^{S_\beta}_{S_\al}}
}
$$
where $\pi^{S_\beta}_{S_\al}$ is a projection and $j^{S_\beta}_{S_\al}$ is a canonical embedding. Then the limit map $\map hX{\cube S}$ (where $S=\bigcup_{\al<\kappa}S_\al$) is a Valdivia embedding.
\end{lm}

\begin{pf} The limit of embeddings is an embedding, so it remains to show that $\Sigma(S)\cap \img hX$ is dense in $\img hX$. Fix a standard basic open set $v\subs \cube S$ such that $\img hX\cap v\nnempty$. Find $\al<\kappa$ such that $\img{h_\al}{X_\al}\cap v\nnempty$ (identifying $\cube{S_\al}$ with a subspace of $\cube S$ via $j_{S_\al}^S$). Increasing $\al$ if necessary, we may assume that $v$ depends on (finitely many) coordinates in $S_\al$, i.e. $v=\inv{(\pi^S_{S_\al})}{\img {\pi^S_{S_\al}}v}$.
Since $h_\al$ is a Valdivia embedding, there is $y\in \Sigma(S_\al)\cap v$ such that $y=h_\al(x)$ for some $x\in X_\al$. For every $\beta\goe\al$ we have $j_{S_\al}^{S_\beta}(y)\in v$ and $\suppt(j_{S_\al}^{S_\beta}(y))=\suppt(y)$. By the assumption, $j_{S_\al}^{S_\beta}(y)= h_\beta i_\al^\beta(x)$ and hence $j_{S_\al}^S(y)=h i_\al(x)$. Thus $h i_\al(x)\in v\cap\Sigma(S)$.
\end{pf}

A retraction $\map rXY$ is {\em simple} if 
\begin{enumerate}
	\item for every nonempty open set $u\subs X$ the image $\img ru$ contains a nonempty $G_\delta$ set;
	\item there exists a map $\map gX{\cube{\aleph_0}}$ such that the diagonal map $\map{r\diag g}X{Y\times\cube{\aleph_0}}$ is one-to-one (i.e. $r$ has a {\em metrizable kernel\/}, see \cite{Haydon}).
\end{enumerate}

Observe that every retraction between metrizable spaces is simple. Also, every open retraction with a metrizable kernel is simple. Fix a space $X$ with a point $p\in X$ which is not $G_\delta$ and define $\map r{X\oplus1}X$ by setting $r(x)=x$ for $x\in X$ and $r(*)=p$, where $*$ is the ``new" isolated point in $X\oplus1$. Then $r$ is a retraction with a metrizable kernel, but it is not simple. Below we prove the announced preservation property for simple retractions.

\begin{tw}\label{szop} Assume $X=\liminv\S$, where $\S=\invsys Xr\kappa\al\beta$ is
a continuous inverse sequence of compact spaces such that $X_0$ is Valdivia compact and each $r^{\al+1}_\al$ is a simple retraction. Then $X$ is Valdivia compact.
\end{tw}

\begin{pf} Let $\rinv i\kappa\al\beta$ be a right inverse for $\S$ (Lemma \ref{prawy}). We shall construct inductively Valdivia embeddings $\map{h_\al}{X_\al}{\cube{S_\al}}$ so that the assumptions of Lemma \ref{smok} are satisfied.
We start with any Valdivia embedding $\map{h_0}{X_0}{\cube{S_0}}$.

Fix $\beta<\kappa$ and assume that $h_\al$ has been defined for each $\al<\beta$. If $\beta$ is a limit ordinal then we can use Lemma \ref{smok}. Assume $\beta=\al+1$. Let $r=r^{\al+1}_\al$, $i=i^{\al+1}_\al$. Note that $x\in \img i{X_\al}$ iff $x=ir(x)$.

Fix a countable set $T$ such that $S_\al\cap T=\emptyset$ and fix $\map g{X_{\al+1}}{\cube T}$ such that $g\diag r$ is one-to-one. We may assume that $g(x)\ne0$ for every $x\in X_{\al+1}$.
Define $\map \phi{X_{\al+1}}{[0,1]}$ by setting
$$\phi(x)=d\Bigl(g(x),(gir)(x)\Bigr),$$
where $d$ is a fixed metric on $\cube T$ which is bounded by $1$.
Observe that $\phi^{-1}(0)=\img i{X_\al}$. Indeed, if $x\notin \img i{X_\al}$ then $i(r(x))\ne x$ and since $r(i(r(x)))=r(x)$, it must be that $g(i(r(x)))\ne g(x)$ which implies $\phi(x)>0$. Let $\Qyu_0=(0,1]\cap \Qyu$. For each $q\in\Qyu_0$ define $\map{\phi_q}{X_{\al+1}}{[0,1]}$ by setting $\phi_q(x)=\max\{0,\phi(x)-q\}$. Finally, set $S_{\al+1}=S_\al\cup(T\times\Qyu_0)$ and define $\map{h_{\al+1}}{X_{\al+1}}{\cube{S_{\al+1}}}$ by setting
$$h_{\al+1}=(h_\al r)\diag f,\qquad\text{ where }f=\Diag_{q\in\Qyu_0}{(\phi_q\cdot g)}.$$ 
We claim that $h_{\al+1}$ is a Valdivia embedding. Fix $x,x'\in X_{\al+1}$. If $r(x)\ne r(x')$ then $h_\al r(x)\ne h_\al r(x')$. If $r(x)=r(x')$ then $g(x)\ne g(x')$ and $\dn x{x'}\not\subs\img i{X_\al}$, because otherwise $x=ir(x)=ir(x')=x'$.
Assuming $\phi(x)\loe q<\phi(x')$ for some $q\in\Qyu_0$, we have $\phi_q(x)\cdot g(x)=0$ and $\phi_q(x')\cdot g(x')\ne0$. Thus in both cases we have $h_{\al+1}(x)\ne h_{\al+1}(x')$ which shows that $h_{\al+1}$ is an embedding.

Fix an open set $v$ such that $\img{h_{\al+1}}{X_{\al+1}}\cap v\nnempty$. Let $u=\inv{h_{\al+1}}v$. Since $\img ru$ contains a nonempty $G_\delta$ set and the set
$$M=\setof{x\in X_\al}{h_\al(x)\in\Sigma(S_\al)}$$ is dense and countably closed, there exists $x\in M\cap \img ru$. Let $x=r(y)$, $y\in u$. Then $h_{\al+1}(y)\in v\cap\Sigma(S_{\al+1})$, because $S_{\al+1}\setminus S_\al=T\times\Qyu_0$ is countable. It follows that $h_{\al+1}$ is a Valdivia map.

By the definition of $h_{\al+1}$ we have
$$\pi^{S_{\al+1}}_{S_\al} h_{\al+1} = h_\al r^{\al+1}_\al\qquad\text{and}\qquad h_{\al+1}i^{\al+1}_\al = (h_\al ri)\diag(fi) = h_\al\diag 0_{T\times\Qyu_0} = j^{S_{\al+1}}_{S_\al} h_\al,$$
where $0_{T\times\Qyu_0}$ denotes the constant $0$ function in $\cube{T\times\Qyu_0}$. Finally, for $\xi<\al$ we have
$$\pi^{S_{\al+1}}_{S_\xi} h_{\al+1} = \pi^{S_\al}_{S_\xi}\pi^{S_{\al+1}}_{S_\al}h_{\al+1} = \pi^{S_\al}_{S_\xi}h_\al r^{\al+1}_\al = h_\xi r^\al_\xi r^{\al+1}_\al = h_\xi r^{\al+1}_\xi$$
and similarly $h_{\al+1}i^{\al+1}_\xi = j^{S_{\al+1}}_{S_\xi}h_\xi$. 
By Lemma \ref{smok}, this completes the proof.
\end{pf}

\begin{wn}\label{maly} Assume $X=\liminv\S$, where $\S=\invsys Xr{\omega_1}\al\beta$ is a continuous inverse sequence such that each $X_\al$ is a compact metric space and each $r^{\al+1}_\al$ is a retraction. Then $X$ is Valdivia compact. \end{wn}

Corollary \ref{maly} implies that a compact space of weight $\loe\aleph_1$ is Valdivia compact iff it is representable as the limit of an inverse sequence of metric spaces whose bonding mappings are retractions. 

Below we give an answer to the questions on retracts and open images, in case where the image has weight $\aleph_1$.

\begin{tw}\label{sdgfaq3rwrrff}
Let $X$ be a Valdivia compact and assume $\map fXY$ is a retraction or $Y$ is 0-dimensional and $f$ is an open surjection. If\/ $\w(Y)\loe\aleph_1$ then $Y$ is Valdivia compact.
\end{tw}

\begin{pf} We assume that $\w(Y)=\aleph_1$. Assuming $X\subs\cube\lam$ and the inclusion is a Valdivia embedding, by Lemma \ref{wiewior1}, we get $S\subs \lam$ such that $|S|=\aleph_1$ and $f(x)=f(x\mid S)$ for every $x\in X$. Thus, taking $X\mid S$ instead of $X$, we may assume that $\w(X)=\aleph_1$. Let $\S=\invsys Xr{\omega_1}\al\beta$ be a continuous inverse sequence of retractions such that $X=\liminv\S$ and each $X_\al$ is second countable. Let $\S'=\invsys Yp{\omega_1}\al\beta$ be the inverse sequence induced by $\S$ and $f$, i.e. $Y=\liminv\S'$ and there exist mappings $\map{f_\al}{X_\al}{Y_\al}$ such that $p_\al f=f_\al r_\al$ holds for every $\al<\omega_1$. 

Assume first that $f$ is a retraction.
By the Spectral Theorem \ref{asfjsifjioqwrqr}, there exists a closed unbounded set $C\subs\omega_1$ such that $f_\al$ is a retraction for every $\al\in C$. If $\al\in C$ then $p_\al f$ is a retraction, which implies that $p_\al$ is a retraction as well. Thus $Y=\liminv\S''$, where $\S''=\S'\rest C$ is an inverse sequence of metric compact spaces with retractions. By Corollary \ref{maly}, $Y$ is Valdivia compact.

Assume now that $f$ is an open surjection and $Y$ is $0$-dimensional. Then, by Proposition \ref{wsegasgpajspjwqprzmv}, there exists a closed unbounded set $C\subs\omega_1$ such that for every $\al\in C$, $Y_\al$ is $0$-dimensional and $f_\al$ is open. Now, if $\al\in C$ then $f_\al$ is an open surjection from a compact metrizable space onto a $0$-dimensional space. Thus $f_\al$ is a retraction and we can repeat the above argument.
\end{pf}

Recall that, according to \cite{KU}, there exists a Valdivia compact space, namely a product of $\aleph_1$ many metric compact groups, which has an open map (in fact: an epimorphism) onto a compact group which is not Valdivia. Thus the assumption on $0$-dimensionality is essential in Theorem \ref{sdgfaq3rwrrff}.

By a {\em (generalized) cube} we mean an arbitrary product of metric compacta. It has been proved in \cite{KL} that retracts of cubes belong to the class of semi-Eberlein compacta. A compact space $K$ is {\em semi-Eberlein} \cite{KL} if $K\subs[0,1]^\kappa$ so that $K\cap c_0(\kappa)$ is dense in $K$, where $c_0(\kappa)$ denotes the well known Banach space of all $\kappa$-sequences convergent to zero. Clearly, every semi-Eberlein space is Valdivia compact. The mentioned result from \cite{KL} is obtained by applying a preservation theorem for semi-Eberlein spaces, similar to Theorem \ref{szop}, with stronger assumptions on the successor bonding maps. 

Our Theorem \ref{szop} gives the weaker statement, namely that 

\begin{wn} Retracts of cubes are Valdivia compact. \end{wn}

\begin{pf} Let $X$ be a retract of a generalized cube. By the result of Shchepin \cite{Szczepin} (or rather by the methods from \cite{Szczepin}), $X=\liminv\S$, where $\S=\invsys Xp\kappa\al\beta$ is a continuous inverse sequence such that $X_0$ is a compact metric space, each $p^{\al+1}_\al$ is an open retraction with a metrizable kernel. 
For the details we refer to \cite{K}.
Thus, each $p^{\al+1}_\al$ is a simple retraction, so $X$ is Valdivia by Theorem \ref{szop}.
\end{pf}

We finish this section by showing that some obvious modifications of the assumptions in Theorem \ref{szop} lead to false statements.

\begin{ex}(a) Let $\S=\invsys Xr{\aleph_2}\al\beta$ be such that $X_\al=\al+1$ and $\map{r^\beta_\al}{\beta+1}{\al+1}$ is a retraction which maps every $\xi>\al$ to $\al$. Then $\S$ is a continuous inverse sequence, each $r^{\al+1}_\al$ is a two-to-one retraction with a metrizable kernel and each $X_\al$ is Valdivia compact.
On the other hand, it is easy to see (e.g. using Lemma \ref{antybrzim}) that $\liminv\S$ is homeomorphic to the linearly ordered space $\omega_2+1$, which is not a continuous image of any Valdivia compact (see \cite{Kalenda6}).

(b) Fix a countable dense set $D=\setof{d_n}{\ntr}$ in $2^{\omega_1}$ and define
$$X=(2^{\omega_1}\times\sn{0_\omega})\cup\setof{e_n}{\ntr},$$
where $0_\omega$ denotes the constant zero function on $\omega$ and $e_n=d_n\concat\chi_{\sn n}$ ($\chi_{\sn n}\in2^\omega$ denotes the characteristic function of $\sn n$). It is easy to see that $X\subs 2^{\omega_1+\omega}$ is homeomorphic to the Alexandrov duplicate of $D$ in $2^{\omega_1}$, i.e.
$$X\iso (2^{\omega_1}\times\sn0)\cup (D\times\sn1),$$
where $D\times \sn1$ consists of isolated points and a basic neighborhood of $\pair p0\in X$ is of the form $(V\times\dn01)\setminus F$, where $V$ is a neighborhood of $p$ in $2^{\omega_1}$ and $F\subs D\times\sn1$ is finite. Thus, $X$ is a non-metrizable compactification of the natural numbers, therefore it is not a continuous image of any Valdivia compact. On the other hand, $X=\liminv\invsys Xr{\aleph_0} nm$, where $X_n=X\mid(\omega_1+n)$ and $r^m_n(x)=x\mid(\omega_1+n)$ for $x\in X_m$ ($n<m$). Each $X_n$ is Valdivia compact (being the union of $2^{\omega_1}$ and a finite set) and each $r^{n+1}_n$ is a retraction with a metrizable kernel (in fact $X_{n+1}\setminus X_n=\sn{e_n}$ and $r^{n+1}_n(e_n)=d_n$).
\end{ex}

\section{Inverse systems with right inverses II}

In this section we investigate $\sig$-complete inverse systems whose bonding mappings are retractions. The aim is to obtain the characterization of Valdivia compacta mentioned in the introduction.

An {\em r-skeleton} in a space $X$ is a pair $\pair\S\II$, where $\S=\invsys Xr\Sigma st$ is a $\sig$-complete inverse system of second countable spaces such that $X=\liminv\S$ and $\II=\rinv i\Sigma st$ is a right inverse of $\S$.
An r-skeleton $\pair\S\II$ is {\em commutative} if
$$(i_s r_s) (i_t r_t) = (i_t r_t) (i_s r_s)$$
holds \emph{for every} $s,t\in \Sigma$. Note that, regardless of commutativity, we always have $(i_s r_s) (i_t r_t) = i_s r_s = (i_t r_t) (i_s r_s)$, whenever $s<t$. 

Let $\pair\S\II$ be an r-skeleton, where $\S$, $\II$ are as above. Define $R_s=i_s r_s$. By Lemma \ref{brzim}, $\sett{R_s}{s\in\Sigma}$ has the following properties:
\begin{enumerate}
	\item[(1)] $s\loe t\implies R_s R_t = R_s = R_t R_s$.
	\item[(2)] $x=\lim_{s\in\Sigma}R_s(x)$ for every $x\in X$.
	\item[(3)] If $T\subs\Sigma$ is countable and directed then $t=\sup T$ exists and $R_t(x)=\lim_{s\in T}R_s(x)$ for every $x\in X$.
\end{enumerate}
Conversely, assume that $\Sigma$ is a directed poset and $\setof{R_s}{s\in\Sigma}$ is a collection of selfmaps of $X$ satisfying conditions (1)--(3) and such that $\img{R_s}X$ is second countable for each $s\in\Sigma$. Setting $X_s=\img{R_s}X$ and $r^t_s=R_t\rest X_t$, we get a $\sig$-complete inverse system $\S=\invsys Xr\Sigma st$ such that $X=\liminv \S$ (see Lemma \ref{antybrzim}) and the inclusions provide a right inverse to $\S$. Thus $\setof{R_s}{s\in\Sigma}$ determines an r-skeleton on $X$.
We shall say that $\setof{R_s}{s\in\Sigma}$ is an {\em internal r-skeleton} in $X$.

We are going to prove that a compact space is Valdivia if and only if it has a commutative r-skeleton. For this aim we need some auxiliary results concerning r-skeletons.

\begin{lm}\label{sosna} Let $\setof{R_s}{s\in\Sigma}$ be an internal r-skeleton in a compact space $X$, let $\Sigma_0\subs\Sigma$ be $\sig$-closed and let $D=\bigcup_{s\in\Sigma_0}X_s$. Then for every map $\map fDY$ into a second countable space $Y$ there exists $t\in\Sigma_0$ such that $f = f R_t$. In particular $\cl_XD=\beta D$.
\end{lm}

\begin{pf}
Fix $\chi$ big enough and a countable elementary substructure $M$ of $H(\chi)$ such that $f,\Sigma_0\in M$ and $\setof{R_s}{s\in\Sigma}\in M$. Let $t=\sup(\Sigma_0\cap M)$. Then $t\in\Sigma_0$. We claim that $f(x)=f(R_t(x))$ for every $x\in D$. Suppose otherwise and fix in $M$ disjoint basic open sets $v,w\subs Y$ such that $f(x)\in v$ and $f(R_t(x))\in w$ for some $x\in D$.
Since $X_t=\liminv\invsys XR{\Sigma_0\cap M}s{s'}$, where $R^{s'}_s=R_s\rest X_{s'}$, there are $s\in\Sigma_0\cap M$ and an open set $u\subs X_s$ such that
\begin{equation}
R_t(x)\in \inv{R_s}u\quad\text{ and }\quad X_t\cap\inv{R_s}u\subs \inv fw.
\tag{*}\end{equation}
In particular $x\in \inv{R_s}u$, because $s<t$ and $R_sR_t(x)=R_s(x)$. We may assume that $u\in M$, because $X_s$ is second countable.
By elementarity, $M\models D\cap \inv{R_s}u\cap \inv fv\nnempty$, because $x\in D\cap \inv{R_s}u\cap \inv fv$. Fix $y\in M\cap D\cap \inv{R_s}u\cap \inv fv$. Then $y\in X_{s'}$ for some $s'\in\Sigma_0\cap M$ and therefore $y\in X_t\cap\inv{R_s}u$, which by (*) implies $f(y)\in w$. Hence $v\cap w\nnempty$, a contradiction. \end{pf}

\begin{lm}\label{sosna2} Assume $\setof{R_s}{s\in\Sigma}$ is an internal r-skeleton in a compact space $X$ and let $T\subs\Sigma$ be a nonempty $\sig$-closed set. Define $\map{R_T}XX$ by setting $R_T(x)=\lim_{t\in T}R_t(x)$. Then $R_T$ is a well defined retraction onto $X_T=\cl_X(\bigcup_{t\in T}X_t)$ and $R_t R_T=R_t$ for every $t\in T$. Moreover $\pair{X_T}{\sett{R_t\rest X_T}{t\in T}}=\liminv\invsys XRT t{t'}$, where $R_t^{t'}=R_t\rest X_T$.
\end{lm}

\begin{pf} Fix $x\in X_T$ and its neighborhood $v$ in $X_T$. Choose $\map f{X_T}{\cube{}}$ such that $f(x)=1$ and $\inv f{(0,1]}\subs v$. By Lemma \ref{sosna}, $f\rest D = (f\rest D)R_t$ for some $t\in T$, where $D=\bigcup_{t\in T}X_t$. Since $D$ is dense in $X_T$, actually $f=fR_t$. Thus, given $t'\goe t$, $t'\in T$, we have $$fR_{t'}(x)=fR_tR_{t'}(x)=fR_t(x)=f(x).$$
Hence $R_{t'}(x)\in v$ for every $t'\goe t$.
Since $v$ was arbitrary, this shows that $x=\lim_{t\in T}R_t(x)$ for every $x\in X_T$. By Lemma \ref{antybrzim}, $\pair{X_T}{\sett{R_t\rest X_T}{t\in T}}=\liminv\invsys XRT t{t'}$.
Now since for every $t<t'$ in $T$ we have $R_t=R^{t'}_t R_{t'}$, there exists a unique map $\map gX{X_T}$ such that $(R_t\rest X_T) g=R_t$ holds for every $t\in T$. Setting $R_T=i g$, where $i$ denotes the inclusion of $X_T$ into $X$, we have
$R_t R_T=R_T$ for every $t\in T$ and consequently
$$R_T(x)=\lim_{t\in T}R_t(R_T(x))=\lim_{t\in T}R_t(x)$$ for every $x\in X$. This also implies that $R_T R_T=R_T$, i.e. $R_T$ is a retraction.
\end{pf}

We shall need the following simple property of $\sig$-complete inverse systems.

\begin{lm}\label{pjietpg}
Assume $\S=\invsys Xp\Sigma st$ is a $\sig$-complete inverse system of compact spaces and $T\subs\Sigma$ is a directed set such that $\Sigma$ is the only $\sig$-closed subset of $\Sigma$ which contains $T$. Then $\liminv\S=\liminv(\S\rest T)$.
\end{lm}

\begin{pf} Fix $x_0,x_1$ in $X=\liminv\S$ and assume $p_t(x_0)=p_t(x_1)$ for every $t\in T$. We have to show that $x_0=x_1$. Define $\Sigma'=\setof{s\in\Sigma}{p_s(x_0)=p_s(x_1)}$. Then $T\subs\Sigma'$ and by the $\sig$-continuity of the system, $\Sigma'$ is $\sig$-closed. Thus $\Sigma'=\Sigma$ and consequently $x_0=x_1$.
\end{pf}

\section{A characterization of the class of Valdivia compacta}

In this section we prove the announced categorical characterization of the class of Valdivia compact spaces:

\begin{tw}\label{w13} A compact space is Valdivia compact if and only if it has a commutative r-skeleton. \end{tw}

Let $F$ be a covariant functor on the category of compact spaces (i.e. $F(K)$ is compact whenever $K$ is compact and $\map{F(f)}{F(K)}{F(L)}$ whenever $\map fKL$).
We say that $F$ is {\em continuous} if it preserves inverse limits and we say that $F$ is {\em weight preserving} if $\w(F(K))\loe \w(K)$ for every compact $K$. Typical examples of continuous weight preserving functors are probability measures, Vietoris hyperspace, symmetric $n$-power (i.e. the hyperspace of at most $n$-element sets) and superextension.

Since the notion of a commutative r-skeleton is defined in the language of category theory, Theorem \ref{w13} implies the following

\begin{wn} The class of Valdivia compact spaces is stable under continuous
weight preserving covariant functors on compact spaces.
\end{wn}

The proof of Theorem \ref{w13} relies on the following technical lemma.

\begin{lm}\label{grzyb} Assume $X$ has a commutative internal r-skeleton $\S=\setof{R_s}{s\in\Sigma}$ and assume $\Ef$ is a finite (possibly empty) collection of retractions of $X$ which commute with $\S$, i.e. $f R_s=R_s f$ for every $f\in\Ef$ and $s\in \Sigma$.
Then there exists a collection of sets $\Yu$ such that, letting $X_s:=\img{R_s}X$, we have
\begin{enumerate}
	\item[(1)] Each element of $\Yu$ is of the form $\inv{R_s}V$ for some $s\in \Sigma$ and an open set $V\subs X_s$.
	\item[(2)] $\Yu$ is $T_0$ separating on $X\setminus \bigcup_{f\in\Ef}\img fX$ and\/ $\bigcup\Yu=X\setminus \bigcup_{f\in\Ef}\img fX$.
	\item[(3)] Each $X_s$ intersects only countably many elements of\/ $\Yu$.	
\end{enumerate}
\end{lm}

\begin{pf} Note that, since each $X_s$ is second countable (and therefore separable), condition (3) is equivalent to saying that every point of $\bigcup_{s\in\Sigma}X_s$ is covered by only countably many elements of $\Yu$.

We use induction on $\kappa=\w(X)$. The statement holds trivially if $\kappa=\aleph_0$. Assume $\kappa>\aleph_0$ and the statement holds for spaces of weight $<\kappa$. Let $Y=\bigcup_{f\in\Ef}\img fX$.

Recall that $X=\liminv\invsys XR\Sigma s{s'}$, where $R^{s'}_s=R_s\rest X_{s'}$. 
Let $\Bee=\sett{V_\al}{\al<\kappa}$ be a base of $X$ such that $V_\al=\inv{R_{s_\al}}{U_\al}$, where $U_\al\subs X_{s(\al)}$ is open ($\al<\kappa$).
It is straight to construct a continuous chain $\sett{T_\al}{\al<\kappa}$ of directed subsets of $\Sigma$ such that $\setof{s_\xi}{\xi<\al}\subs T_\al$ and $|T_\al|<\kappa$ for every $\al<\kappa$. Let $T=\bigcup_{\al<\kappa}T_\al$. Then $T$ is a directed subset of $\Sigma$. Let $\map pX{X_T}$ be the canonical projection, where $X_T=\liminv\invsys XRTs{s'}$. Then $p$ is one-to-one, because $\setof{s_\al}{\al<\kappa}\subs T$. Thus $X=X_T$ and, by Lemma \ref{pjietpg}, we may assume that $\Sigma$ is the smallest $\sig$-closed set containing $T$.
Let $\Sigma_\al\subs \Sigma$ be the smallest $\sig$-closed set which contains $T_\al$. Then $\sett{\Sigma_\al}{\al<\kappa}$ is a chain of $\sig$-closed subsets of $\Sigma$ such that
\begin{enumerate}
	\item[(i)] $X_\al:=\liminv\invsys XR{\Sig_\al}s{s'}$ has weight $<\kappa$ for each $\al<\kappa$.
	\item[(ii)] If $\delta<\kappa$ is a limit ordinal then $\Sig_\delta$ is the smallest $\sig$-closed set containing $\bigcup_{\al<\delta}\Sig_\al$.
  \item[(iii)] $\Sigma$ is the smallest $\sig$-closed set containing $\bigcup_{\al<\kappa}\Sigma_\al$.
\end{enumerate}
Property (i) follows from Lemma \ref{pjietpg} and the fact that $|T_\al|<\kappa$. Property (ii) follows from the continuity of the chain $\sett{T_\al}{\al<\kappa}$. 

By Lemma \ref{sosna2}, we may assume that $X_\al=\cl(\bigcup_{s\in\Sig_\al}X_s)$ and the projection $\map{R_\al}X{X_\al}$ satisfies $R_\al(x)=\lim_{s\in\Sig_\al}R_s(x)$ for every $x\in X$ (to avoid confusion, we assume that $\Sig$ does not contain ordinals). Given $t\in\Sigma$ we have
$$R_t R_\al(x)= R_t(\lim_{s\in\Sig_\al} R_s(x)) = \lim_{s\in\Sig_\al} R_tR_s(x) = \lim_{s\in\Sig_\al} R_sR_t(x) = R_\al R_t(x).$$
Thus
\begin{enumerate}
	\item[(iv)] $R_\al R_t = R_t R_\al$ for every $\al<\kappa$, $t\in\Sigma$.
\end{enumerate}
Similarly, $R_\al R_\beta=R_\beta R_\al=R_\al$ whenever $\al<\beta<\kappa$. Now (ii) and (iii) together with Lemma \ref{pjietpg} imply that
\begin{enumerate}
	\item[(v)] $X=\liminv\invsys XR\kappa{\al}{\al'}$ and $X_\delta=\liminv\invsys XR\delta{\al}{\al'}$ for a limit ordinal $\delta<\kappa$, where $R^{\al'}_\al=R_\al\rest X_{\al'}$.
\end{enumerate}

We construct inductively an increasing sequence $\setof{\Yu_\al}{\al\loe\kappa}$ of families of subsets of $X$ satisfying the following conditions (where we set $X_\kappa=X$):
\begin{enumerate}
	\item[(a)] Each element of $\Yu_\al$ is of the form $\inv{R_s}V$ for some $s\in\Sig_\al$ and an open set $V\subs X_s$.
	\item[(b)] $\Yu_\al$ is $T_0$ separating on $X_\al\setminus Y$ and $X_\al\cap\bigcup\Yu_\al=X_\al\setminus Y$.
	\item[(c)] For each $s\in\Sig$, $X_s$ intersects only countably many elements of $\Yu_\al$.
	\item[(d)] $U\cap X_\al=\emptyset$ whenever $U\in\Yu_\beta\setminus \Yu_\al$.
\end{enumerate}
Concerning condition (a), we shall use the following observation:
\begin{enumerate}
	\item[(vi)] If $V\subs X_s$, where $s\in\Sigma$, then $V\cap Y=\emptyset$ iff $\inv{R_s}V\cap Y=\emptyset$.
\end{enumerate}
To see (vi), suppose $y\in \inv{R_s}V\cap \img fX$, where $f\in \Ef$. Then $R_s(y)\in V$ and $f(R_s(y))=R_s(f(y))=R_s(y)$ which implies that $R_s(y)\in Y$, i.e. $R_s(y)\in V\cap \img fX$. Conversely, if $y\in V\cap \img fX$ then $y\in X_s$ and hence $R_s(y)=y$, which gives $y\in\inv{R_s}V\cap \img fX$. This shows (vi). These arguments used only the fact that each $f\in\Ef$ commutes with each $R_s$, $s\in\Sigma$. Thus, by (iv), we also have
\begin{enumerate}
	\item[(vii)] If $V\subs X_s$, $s\in\Sigma$ and $\al<\kappa$, then $V\cap X_\al=\emptyset$ iff $\inv{R_s}V\cap X_\al=\emptyset$.
\end{enumerate}

Fix $\beta\loe\kappa$ and suppose that $\Yu_\al$ have already been constructed for $\al<\beta$. In order to avoid repeating the same arguments twice, we set $X_{-1}=\emptyset=\Yu_{-1}$ and we treat $0$ as a successor ordinal. We consider two cases.

{\it Case 1.~} $\beta=\al+1$. Since $\w(X_{\al+1})<\kappa$, we can apply the inductive hypothesis to $X_{\al+1}$ and $\sett{f\rest X_{\al+1}}{f\in\Ef}\cup\sn{R_\al\rest X_{\al+1}}$, obtaining a family $\Vee$ satisfying (1)--(3). Thus, each element of $\Vee$ is of the form $U=\inv{(R_s\rest X_{\al+1})}W=\inv{R_s}W\cap X_{\al+1}$, where $s\in \Sig_{\al+1}$, $W$ is open in $X_s$ and disjoint from $Y$ (by (vi)). Moreover, $\Vee$ is $T_0$ separating on $X_{\al+1}\setminus(X_\al\cup Y)$ and $\bigcup\Vee=X_{\al+1}\setminus(Y\cup X_\al)$.
Define $U^*=\inv{R_s}W$ and set
$$\Yu_{\al+1}=\Yu_\al\cup\setof{U^*}{U\in \Vee}.$$ 
Clearly, (a) is satisfied. Note that $U^*\cap X_{\al+1}=U$. Thus $X_{\al+1}\cap (\Yu_{\al+1}\setminus\Yu_\al)=X_{\al+1}\setminus(X_\al\cup Y)$ and $\Yu_{\al+1}\setminus\Yu_\al$ is $T_0$ separating on $X_{\al+1}\setminus(X_\al\cup Y)$.
By the inductive hypothesis, $\Yu_{\al}$ is $T_0$ separating on $X_\al\setminus Y$. Hence $\Yu_{\al+1}$ is $T_0$ separating on $X_{\al+1}\setminus Y$, because the family $\Yu_{\al+1}\setminus \Yu_\al$ separates every point of $X_{\al+1}\setminus (X_\al\cup Y)$ from every point of $X_\al\setminus Y$. This shows (b). 

In order to show (c), fix $s\in \Sig$ and fix $x\in X_s$. By Lemma \ref{sosna}, applied to the r-skeleton $\setof{R_s\rest X_{\al+1}}{s\in\Sig_{\al+1}}$, we can find $t\in\Sig_{\al+1}$ such that $R_s\rest X_{\al+1}= (R_s\rest X_t) (R_t \rest X_{\al+1})$. Thus we get
$$R_s R_{\al+1}(x)=R_s R_t R_{\al+1}(x)=R_{\al+1}R_t R_s(x)=R_{\al+1}R_t(x)=R_t(x).$$
The last equality follows from the fact that $R_t(x)\in X_t\subs X_{\al+1}$.
On the other hand, we have $R_s R_{\al+1}(x)=R_{\al+1} R_s(x)=R_{\al+1}(x)$. Hence $R_{\al+1}(x)=R_t(x)\in X_t$.

Now, if $x\in U^*=\inv{R_r}W$, where $U=\inv{(R_r\rest X_{\al+1})}W\in\Vee$ and $r\in \Sig_{\al+1}$ then, since $X_r\subs X_{\al+1}$, we have
$$R_r R_{\al+1}(x)=R_{\al+1} R_r(x) = R_r(x)\in W,$$
which shows that $R_{\al+1}(x)=R_t(x)\in U\cap X_t$. By assumption, the set $\setof{U\in\Vee}{R_{t}(x)\in U}$ is countable. Thus $x$ belongs to countably many elements of $\Yu_{\al+1}\setminus \Yu_\al$. This, together with the inductive hypothesis, shows (c).

Condition (d) follows directly from (vii) and from the fact that $X_\al\cap\bigcup\Vee=\emptyset$.

{\it Case 2.~} $\beta$ is a limit ordinal. 
Define $\Yu_\beta=\bigcup_{\al<\beta}\Yu_\al$. It is necessary to check conditions (b) and (c) only.

We first check (b). Fix $x_0\ne x_1$ in $X_\beta\setminus Y$. By (v), there is $\al<\beta$ such that $R_\al(x_0)\ne R_\al(x_1)$ and by Lemma \ref{brzim} we may assume that $R_\al(x_i)\notin Y$ for $i=0,1$. Find $U\in \Yu_\al$ such that e.g. $R_\al(x_1)\in U$ and $R_\al(x_0)\notin U$. Since $U=\inv{R_s}V$ for some $s\in \Sig_\al$, we have $R_sR_\al=R_s$ and therefore $x_1\in \inv{(R_sR_\al)}V=\inv{R_s}V=U$ and by the same reason $x_0\notin U$.
This shows (b).

Condition (c) holds trivially if $\cf\beta=\aleph_0$, so assume $\cf\beta>\aleph_0$, which implies that $\Sig_\beta=\bigcup_{\al<\beta}\Sig_\al$. Fix $s\in\Sig$ and fix $x\in X_s$. By Lemma \ref{sosna}, there exists $t\in \Sig_\beta$ such that $R_s\rest X_\beta = (R_s\rest X_t) (R_t\rest X_\beta)$. Hence
$$R_s R_\beta(x)=R_s R_t R_\beta(x)=R_\beta R_t R_s(x)=R_\beta R_t(x)=R_t(x),$$
where the last equality follows from the fact that $X_t\subs X_\beta$. On the other hand, $R_s R_\beta(x)=R_\beta R_s(x)=R_\beta(x)$. Thus $R_\beta(x)=R_t(x)$. Let $\al<\beta$ be such that $t\in\Sigma_\al$.
If $x\in U$ and $U=\inv{R_r}V\in \Yu_\xi$, where $r\in\Sigma_\xi$ and $\xi<\beta$, then 
$$R_r R_t(x)=R_r R_\beta(x)=R_\beta R_r(x)=R_r(x)\in V,$$
which means that $R_t(x)\in U\cap X_t\subs U\cap X_\al$.  By the inductive hypothesis, we have $U\cap X_\al=\emptyset$ for $U\in\Yu_\xi\setminus\Yu_\al$. Hence $\xi\loe\al$. It follows that $x\notin U$ whenever $U\in\Yu_\beta\setminus \Yu_\al$. This, together with the inductive hypothesis, shows (c). 

It follows that the construction can be carried out. Setting $\Yu=\Yu_{\kappa}$ we obtain the desired collection.
\end{pf}

\begin{pf}[Proof of Theorem \ref{w13}.]
Assume $X$ has an internal commutative r-skeleton $\setof{R_s}{s\in\Sig}$. Applying Lemma \ref{grzyb} with $\Ef=\emptyset$, we obtain a $T_0$ separating collection $\Yu$ which consists of open $F_\sig$ sets and each point of $D=\bigcup_{s\in\Sig}X_s$ belongs to only countably many elements of $\Yu$.
Taking for each $u\in\Yu$ a function $\map{h_u}X{[0,1]}$ such that $u=\inv{h_u}{(0,1]}$, the diagonal map $\Diag_{u\in\Yu}h_u$ provides a Valdivia embedding of $X$ into $\cube{\Yu}$.

Assume now that $X$ is Valdivia compact and the inclusion $X\subs[0,1]^\kappa$ is a Valdivia embedding. Let $\Sig=\setof{S\in\dpower\kappa{\aleph_0}}{X\mid S\subs X}$ and define $R_S(x)=x\mid S$. By Lemma \ref{wiewior1}, $\setof{R_S}{S\in\Sig}$ is an internal r-skeleton in $X$. It is clear that $R_S R_T=R_T R_S$ holds for every $S,T\in\Sig$.
\end{pf}

The following example shows that the assumption on commutativity in Theorem \ref{w13} cannot be omitted.

\begin{ex} Let $\kappa$ be an infinite cardinal. Denote by $\Sigma$ the collection of all closed countable sets $A\subs\kappa+1$ such that $0\in A$ and every isolated point of $A$ is isolated in $\kappa+1$. Then $\pair{\Sig}{\subs}$ is a directed poset. For each $A\in\Sig$ define $\map{r_A}{\kappa+1}{\kappa+1}$ by setting $r_A(x)=\max\setof{a\in A}{a\loe x}$. The assumption on isolated points implies that $r_A$ is continuous, so it is a retraction onto $A$. Furthermore, if $A,B\in\Sig$ and $A\subs B$ then $r_Ar_B = r_A = r_Br_A$. Finally, $x=\lim_{A\in\Sig}r_A(x)$ and also $r_S(x)=\lim_{A\in\Aaa}r_A(x)$ whenever $S=\sup\Aaa$ and $\Aaa\subs \Sig$ is countable and directed. Thus $\setof{r_A}{A\in\Sig}$ is an internal r-skeleton in $\kappa+1$. On the other hand, if $\kappa>\aleph_1$ then $\kappa+1$ is not Valdivia compact.
\end{ex}

\noindent{\bf Acknowledgments.} The first author would like to thank The Fields Institute for Research in Mathematical Sciences (Toronto, Canada, September -- December 2002) and the University of Prince Edward Island (Charlottetown, Canada, January -- April 2003) for their great hospitality, when the main part of this research was being done.

\end{document}